\documentclass[11pt]{amsart}
\usepackage{amsmath,amsfonts,amssymb}
\let\L=\Lambda
\let\l=\lambda

\let\g=\gamma
\let\t=\tau
\let\S=\Sigma
\let\e=\epsilon
\let\O=\Omega
\let\bb=\mathbb
\let\d = \partial
\let\cl=\mathcal

\let\z=\zeta

\newtheorem{theorem}{Theorem}[section]

\newtheorem{lemma}[theorem]{Lemma}
\newtheorem{corollary}[theorem]{Corollary}
\newtheorem{definition}[theorem]{Definition}
\begin{document}
\title{Holomorphic correspondences between CR manifolds.}
\author{C. Denson Hill and Rasul Shafikov}
\maketitle

\section{Introduction.}

One of the interesting phenomena in several complex variables is the
analytic continuation of a germ of a biholomorphic map $f: M\to M'$,
defined at a point $p\in M$, where $M, M'\subset \bb C^n$ are
real-analytic hypersurfaces. Already Poincar\'e \cite{po} observed
that a biholomorphic map sending an open piece of a unit sphere in
$\bb C^2$ to another such open piece must be an automorphism of the
unit ball. This was proved for $\bb C^n$ in \cite{ta} and
\cite{al}. Clearly such an extension is possible only  for $n>1$, thus
showing the very special nature of CR maps between real-analytic CR
manifolds. 

Pinchuk \cite{pi1},\cite{pi3} proved that the germ of a non-degenerate
holomorphic map from a real-analytic strictly pseudoconvex
hypersurface $M$ in $\bb C^n$ to a unit sphere in $\bb C^N$, $N\ge n$, 
extends as a holomorphic map along any path on $M$. Webster \cite{w}
proved that the germ of a biholomorphic map 
between real-algebraic Levi non-degenerate hypersurfaces in $\bb C^ n$ 
extends as an algebraic map, and also gave sufficient conditions for 
the map to be rational. In that paper he studied the geometric
properties of Segre varieties, which were originally introduced in
\cite{se}.  

Much attention was devoted to the generalization of Webster's theorem
to the case of different dimensions and higher codimensions. In this
situation both $M$ and $M'$ are assumed 
to be real algebraic submanifolds or sets; that is, defined by the
zero locus of a system of real polynomials. Under certain conditions,
it then turns out that a locally defined holomorphic map between such
objects must necessarily be algebraic. See Baouendi, Ebenfelt and
Rothschild \cite{ber}, Huang\cite{hu}, Sharipov and Sukhov
\cite{shsu},  Baouendi, Huang and Rothschild \cite{bhr}, Coupet,
Meylan and Sukhov \cite{cms}, Zaitsev \cite{z}, Merker \cite{me} and
many additional references contained therein. The special case of
hyperquadrics was also considered in Tumanov \cite{tu}, Forstneri\v c
\cite{fo}, Sukhov \cite{su} and other papers.  

On the other hand, much less is known if at least one of the
submanifolds is not assumed to be real-algebraic. In this case the map
need not be algebraic; however, analytic continuation along a
real-analytic hypersurface is also possible. Pinchuk \cite{pi2} proved
that a germ 
of a biholomorphic map $f$ from a strictly pseudoconvex,
real-analytic, non-spherical hypersurface $M$ to a compact, strictly
pseudoconvex, real-analytic hypersurface $M'\subset \bb C^n$ extends
holomorphically along any path on $M$. A similar result was shown in
\cite{sh1} for the case when $M$ is essentially finite, smooth,
real-analytic and $M'\subset \bb C^n$ is compact, real-algebraic and
strictly pseudoconvex. Levi non-degeneracy of the target hypersurface
ensures that the extended map is single valued. If the target
hypersurface is just assumed to be compact and smooth real-algebraic,
the extension in general will be multiple-valued as was proved in
\cite{sh}. This naturally leads to consideration of holomorphic
correspondences, a multiple-valued generalization of holomorphic
mappings. 

In this paper we study the analytic continuation of germs of holomorphic
mappings from smooth real-analytic CR submanifolds of arbitrary
codimension to compact smooth real-algebraic generic submanifolds
in $\bb P^N$ of general codimension. The continuation
that we obtain is a holomorphic correspondence from a
neighborhood of the submanifold in the pre-image to $\bb P^N$. This is
analogous to the algebraicity of the map asserted in the case when both
submanifolds are real-algebraic. We also study some applications to
maps between pseudoconcave CR submanifolds in $\bb P^n$. It is rather
surprising that under certain conditions, a local CR map between such
objects turns out to be the restriction of a rational, or even a linear
map in $\bb P^n$ without the assumption of algebraicity of the
submanifold in the pre-image. 

Our results generalize the extension property of a germ of a 
biholomorphic mapping from a compact real-analytic hypersurface in
$\mathbb C^n$ to a compact real-algebraic hypersurface in 
$\mathbb C^n$ proved in \cite{sh}. We remark that the proofs of 
the main results of this paper differ significantly from those 
utilized in \cite{sh}, where the main construction essentially uses 
the fact that the Segre varieties have codimension one.

In the next section we present the main results. In Section 3 we give
some background on CR manifolds, Segre varieties and holomorphic
correspondences. Section 4 contains the proof of the local extension
of a holomorphic map as a holomorphic correspondence. In Section 5 we
prove the global extension. The last section contains applications of
the main theorem to pseudoconcave CR submanifolds in $\bb P^n$.

The second author would like to acknowledge the support of the
Max-Planck-Institute for mathematics in Bonn during the preparation of
this paper. We would also like to thank S.~Popescu and S.~Lu for a
number of very helpful discussions.

\section{Statement of Results.}

\subsection{Analytic Continuation}
By a holomorphic correspondence we mean a complex analytic subset
$A\subset X\times X'$, where $X$ and $X'$ are complex manifolds, such
that (i) $\dim A \equiv \dim X$, and (ii) the projection $\pi: A\to X$
is proper. It is natural to define a multiple valued map $F=\pi'\circ
\pi^{-1}: X\to X'$ associated with $A$. A holomorphic correspondence
$F$ is called a {\it finite} correspondence, if $F(p)$ and
$F^{-1}(p')$ are finite sets for any $p\in X$ and $p'\in X'$. We say
that $F$ {\it splits} at a point 
$q\in X$ if there exist a neighborhood $U_q$ of $q$ and an integer $k$
such that $F|_{U_q}$ is just the union of $k$ holomorphic maps $F^j:
U_q\to X'$, $j=1,\dots, k$.

Recall that if $M$ is an abstract paracompact real analytic CR
manifold of type $(m,d)$, then by \cite{af} there exists a
complex manifold $X$ of dimension $n=m+d$ such that $M$ can be
generically embedded into $X$ as a CR submanifold. 

Consider the following situation: $M$  is a smooth real-analytic
minimal CR manifold of type $(m,d)$, $m,d>0$; $M'\subset \mathbb P^N$
is a compact smooth real-algebraic essentially finite generic
submanifold of type $(m,d')$, $d'>0$. The main results of the paper
are the following.   

\begin{theorem}\label{t0}
  Suppose that $\omega\subset M$ is a relatively compact connected
  open subset, and $f:\omega\to M'$ is a real-analytic CR map such
  that $df|_{H_p}: H_p M  \to H_{f(p)} M'$  is an isomorphism between
  the holomorphic tangent spaces for almost all $p\in \omega$. Let  
  $q\in \d \omega$. Then there exists a neighborhood $U_q\subset X$ of
  $q$ such that $f$ extends to $U_q$ as a holomorphic correspondence
  $F:U_q\to \bb P^N$ with $F(U_q\cap M)\subset M'$.
\end{theorem}

We remark that in Theorem~\ref{t0} we may assume that $f:\omega\to M'$
is a smooth CR mapping and the characteristic variety $\nu_z(f)$ has
dimension zero for $z\in \omega$, since by the result in \cite{d}, $f$ is
in fact real-analytic. Also note that we do not claim that the extended 
correspondence near $q$ is finite. 

The Segre map associated with a real-analytic CR manifold $M$ is defined 
by $\lambda :w \to Q_w$, where $Q_w$ is the Segre variety of a point $w$.
We say that $\lambda$ is locally injective at a point $q\in M$, if 
there exists a small neighborhood $U_q$ of $q$ such that $\lambda$ is 
an injective map from $U_q$ onto its image. For details see Section 3. 

\begin{theorem}\label{t1} Assume in addition that $M$ is
  essentially finite and $d'\ge d$. Let $p\in M$, $U_p\subset X$ be a
  neighborhood of $p$, and let $f:U_p\to \bb P^N$ be a holomorphic
  mapping of maximal rank such that $f(U_p\cap M)\subset M'$. Suppose
  that  $M_1$ is a relatively compact simply-connected open subset of
  $M$ containing $p$. Then
  \begin{enumerate}
    \item There exists a neighborhood $U\subset X$ of $M_1$ such that
      $f$ extends as a finite holomorphic correspondence $F: U\to \bb
      P^N$ with $F(M_1)\subset M'$. 
    \item $F$ splits into holomorphic mappings at
    every point $q\in M_1\setminus F^{-1}(\Sigma')$, where $\Sigma'$
    is the set of points on $M'$ near which the Segre map is not
    locally injective. 
  \end{enumerate}
     
\end{theorem}

We note that the assumption in Theorem \ref{t1} that the map $f$ is of 
maximal rank ensures that the images of Segre varieties $Q_z$ under 
$f$ and subsequent analytic continuations of $f$ have the same dimension
as $Q'_{z'}$. The essential finiteness of $M'$ then guaranties that the 
constructed extensions of the map $f$ are finitely-valued. 

To illustrate the conclusions of Theorem \ref{t1} we consider a
simple example. Let $r$ and $s$ be coprime positive integers, and $M$
and $M'$ be the CR hypersurfaces in $\bb P^n$ given in homogeneous
coordinates  by 
\begin{equation}
  M=\{z\in \bb P^n : |z_0|^{2r} + \dots + |z_k|^{2r} -
  |z_{k+1}|^{2r}-\dots -|z_n|^{2r}=0\},
\end{equation}
\begin{equation}
  M'=\{z'\in \bb P^n : |z'_0|^{2s} + \dots + |z'_k|^{2s} -
  |z'_{k+1}|^{2s}-\dots -|z'_n|^{2s}=0\},
\end{equation}
with $k<n$. The finite holomorphic correspondence $F(z)=(z_0^{r/s},
\dots,z_n^{r/s})$ maps $M$ to $M'$, is $r^n:s^n$ valued, and $F$
splits into $s^n$ holomorphic mappings outside the branching locus of
$F$. If we choose a point $p\in M$ outside the branching locus, and if
at $p$ we choose the germ of one of the $s^n$ branches of $F$, then
Theorem \ref{t1} reproduces the entire correspondence $F$.

Theorems \ref{t0} and \ref{t1} can be considered as a
generalization of the results on algebraicity of a local CR map
between real-algebraic submanifolds. Non-algebraic holomorphic
correspondences can be easily constructed from the examples considered
in \cite{bs}. 

Our proof of Theorems \ref{t0} and \ref{t1} is based on the technique
of Segre varieties. As it was mentioned in the introduction, Webster
was the the first to use Segre varieties in the context of holomorphic
mappings. His ideas were further developed in \cite{dw}, \cite{df2},
\cite{dfy}, \cite{dp1} and other papers. Our main construction is
based on the technique of these papers. A somewhat different approach
was developed in \cite{bjt}, \cite{br} and \cite{ber}. We rely on the
characterization of minimality in terms of Segre sets proved in
\cite{ber}. 

Let us briefly explain the idea of the main construction. Let
$p=0\in M$, and let $U_0$ be a small neighborhood of the origin, where 
the map $f$ is defined. The first step is to show that $f$ extends as
a holomorphic correspondence $F_1$ to a neighborhood $U_1$ of $Q_0$,
the Segre variety of the origin. This can be understood as follows. If
$w\in Q_0$, then $0\in Q_w$, and $f(U_0\cap Q_w)$ is a complex
subvariety in the target space which passes through $f(0)$. If
$f(U_0\cap Q_w)\subset Q'_{w'}$ for some $w'$, then we set
$F_1(w)=w'$. For $w$ close to the origin, $f(Q_w\cap U_0)\subset
Q'_{f(w)}$ by the invariance property of Segre varieties, and thus the
extension defined in this way agrees with $f$ on $U_0$. For essentially
finite manifolds, if $w'$ is sufficiently close to $M'$, then there
exists a finite number of points which have the same Segre variety as
$w'$. Therefore in general $F_1$ may not be single valued. The
graph of $F_1$ can be described as  
\begin{equation}
  A_1=\{(w,w')\in U_1\times \bb P^N : f(Q_w\cap U_0)\subset Q'_{w'}\}.
\end{equation}

The next step of the proof is to define inductively the analytic
continuation of $f$ to bigger sets. This can be achieved using Segre
sets (see the next section for definitions). Briefly, once the
analytic continuation of $f$ as a correspondence $F_1$ is established,
we can use a similar construction to define the extension of $F_1$
along the Segre varieties of points on $Q_0$. Let
\begin{equation}\label{prelim2}
  A_2=\{(w,w')\in U_2\times \bb P^N: F_1(Q_w\cap U_1)\subset
  Q'_{w'}\}, 
\end{equation}
where $U_2$ is a neighborhood of $Q^2_0=\cup_{w\in Q_0} Q_w$.
We note that for $w'$ close to $M'$ the condition for $(w,w')$ in 
\eqref{prelim2}  roughly speaking means that all
branches of $F_1$ map $Q_w$ into $Q'_{w'}$. We can repeat
this process for $A_3$, using the extension $A_2$, etc.

Finally, the crucial observation is that after a certain number of
iterations, we obtain an extension of $f$ to a neighborhood of the
origin, which is independent of the choice of $U_0$, where $f$ was
originally defined. This follows from the minimality of $M$ and the
property of Segre sets proved in \cite{ber}. Furthermore, for compact
minimal manifolds such neighborhoods can be chosen of uniform size for
all points on $M$, which allows us to extend $f$ to any
simply-connected relatively compact subset of $M$.

The question remains open whether a similar continuation of $f$ is
possible in the case when $M'$ is real-analytic. This is unknown even 
in the case when $M$ is a pseudoconvex hypersurface and $M'$ is a 
strictly pseudoconvex hypersurface in $\mathbb C^n$. Our method heavily 
relies on the fact that the Segre varieties associated with $M'$ are 
globally defined in $\mathbb P^n$, and therefore cannot be directly 
extended to the real-analytic case.

\subsection{Pseudoconcave CR manifolds}
Our main application of Theorem \ref{t1} concerns pseudoconcave
CR submanifolds embedded into projective spaces. We recall that a CR
manifold is called pseudoconcave, if at each point its Levi form has
at least one positive and at least one negative eigenvalue in every
characteristic conormal direction. 

\begin{theorem}\label{t1.5} Consider a connected $C^{\infty}$-smooth
  compact pseudoconcave CR submanifold $M$ of $\bb P^n$, having type
  $(m,d)$ with $m,d>0$.
  \begin{enumerate}
    \item[(a)] Let $f: M\to \bb P^N$ be a continuous CR map. Then $f$
    is the restriction to $M$ of a rational map $F:\bb P^n\to \bb
    P^N$.  
    \item[(b)] Assume $n=m+d$, so that $M$ is generic in $\bb
    P^n$. Let $f:M\to \bb P^n$ be a CR map which is a local
    diffeomorphism onto $f(M)$. Then $f$ is the restriction of a
    linear automorphism of $\bb P^n$.  
  \end{enumerate}
\end{theorem}

For a locally biholomorphic map from a compact smooth pseudoconcave
hypersurface in $\bb P^n$ to $\bb P^n$ part (b) was first shown in
\cite{iv}. Our proof of Theorem \ref{t1.5} is based on the result of
\cite{hn2}, where it is shown that any CR meromorphic function on $M$
is necessarily rational.    

Combining Theorem \ref{t1.5} with Theorem \ref{t1}, we obtain the
following results. 

\begin{theorem}\label{t2}
  Let $M\subset \bb P^n$ (resp. $M'\subset \bb P^N$) be a
  compact smooth real-analytic pseudoconcave  essentially finite CR 
  submanifold of type $(m,d)$ (resp. $(m,d')$), $m,d,d' > 0$ and
  $d'\ge d$. Let $M$ be simply connected and let $M'\subset \bb P^N$
  be generic real-algebraic, and such that the Segre map associated
  with $M'$ is locally injective. Let $p\in M$, $U_p$ be a
  neighborhood of $p$ in $\bb P^n$, and
  let $f: U_p\cap M\to M'$ be a germ of a smooth CR map of maximal
  rank. Then
  \begin{enumerate}
    \item[(a)] $f$ is the restriction of a rational map $F: \bb P^n \to \bb
    P^N$; 
    \item[(b)] If $n=N=m+d$ and the Segre map associated with $M$ is
    locally injective, then $f$ is the restriction of a linear
    automorphism of $\bb P^n$.
  \end{enumerate}
\end{theorem}

If $M,M'\subset \bb P^n$ are both hypersurfaces, then $f$ may be 
assumed to be just a local CR homeomorphism, since in this case $f$ 
extends smoothly to a neighborhood of $p$, and by \cite{hn1}, the 
Jacobian of $f$ does not vanish at $p$. 

\begin{corollary}\label{c3}
  If a real-analytic submanifold $M$, satisfying the conditions of
  Theorem \ref{t2}, is locally CR equivalent to a real-algebraic CR
  submanifold $M'$ satisfying Theorem \ref{t2}, then $M$ is
  necessarily real-algebraic. 
\end{corollary}

\noindent {\bf Remark.} The above results concerning pseudoconcave CR
manifolds hold also under a weaker assumption on $M$ and $M'$, namely,
instead of pseudoconcavity it is enough to assume that $M$ and $M'$
satisfy the so-called Property E. For a generic $M$ this means that
for each $p\in M$ every local CR function defined near $p$ extends to
a holomorphic function in a full neighborhood of $p$ in the ambient 
space. For details see \cite{hn2}.

Furthermore, in the special case when $M'$ is a hyperquadric we obtain 
the following.

\begin{theorem}\label{c}
  Let $M$ be a simply-connected compact smooth real-analytic Levi
  non-degenerate CR manifold of type $(n-1,1)$, $n>1$. Let $M'$ be the
  hyperquadric in $\bb P^n$ given in homogeneous coordinates by  
  \begin{equation}\label{quadric}
    |z_0|^2+|z_1|^2+\dots+|z_k|^2-|z_{k+1}|^2-\dots-|z_n|^2=0.
  \end{equation}
  Suppose that $\omega$ is a connected open set in $M$, and $f:
  \omega\to M'$ is a CR map that is a local homeomorphism. Then $M$
  and $M'$ are 
  globally CR equivalent; hence $M$ has a CR embedding as a
  hypersurface in $\bb P^n$. In the special case where $0<k<n-1$, and
  $M$ is apriori a hypersurface in $\bb P^n$, then $f$ is the
  restriction to $M$ of a linear automorphism of $\bb P^n$.
\end{theorem} 

Note that Theorem \ref{c} includes both the case when $M'$ is a sphere
in $\bb C^n$, which was proved before in \cite{pi1}, and the case when
$M'$ is a compact pseudoconcave hyperquadric. For the sphere our method
gives an alternative and independent proof of this well known result.

If $M$ is not assumed to
be simply-connected, then $f$ in general may not extend to a global
map from $M$ to $M'$ as examples in \cite{bs} show for the case when
$M'$ is a sphere.

\section{CR manifolds, Segre varieties and holomorphic
  correspondences.} 

An abstract smooth CR manifold of type $(m,d)$ consists of a connected
smooth paracompact manifold $M$ of dimension $2m+d$, a smooth
subbundle $HM$ of $TM$ of rank $2m$, which is called the holomorphic
tangent space of $M$, and a smooth complex structure $J$ on the fibers
of $HM$. Let $T^{0,1}M$ be the complex subbundle of the
complexification $\bb C HM$ of $HM$, which corresponds to the $-i$
eigenspace of $J$:
\begin{equation}
  T^{0,1}M=\{Y+iJY| Y\in HM\}.
\end{equation}
We also require that the formal integrability condition
\begin{equation}
  [C^{\infty}(M,T^{0,1}M),C^{\infty}(M,T^{0,1}M)]\subset
  C^{\infty}(M,T^{0,1}M)  
\end{equation}
holds. We call $m$ the CR dimension of $M$ and $d$ the CR
codimension. $M$ is called minimal at $p\in M$, 
if there exists no local CR manifold $N\subset M$ passing through $p$
having CR dimension $m$, but strictly smaller real dimension. 

The characteristic bundle $H^0M$ is defined to be the
annihilator of $HM$ in $T^*M$. Its purpose is to parametrize the Levi
form, which for every $p\in M$ is defined for $v\in H^0_p M$ and $Y\in
H_p M$ by
\begin{equation}
  \cl L(v; Y)= d\tilde v(Y,JY)=\langle v,[J\tilde Y,\tilde Y]\rangle,
\end{equation}
where $\tilde v\in C^{\infty}(M, H^0 M)$ and $\tilde Y\in (M,HM)$ are
smooth extensions of $v$ and $Y$. For each fixed $v$ it is a Hermitian
quadratic form for the complex structure $J_p$ on $H_p M$. 
A CR manifold $M$ is said to be pseudoconcave if the Levi from $\cl
L(v,\cdot)$ has at least one negative and one positive eigenvalue for
every $p\in M$ and every nonzero $v\in H^0_p M$.

Let $M$ and $M'$ be two abstract smooth CR manifolds, with holomorphic
tangent spaces $HM$ and $HM'$. A smooth map $f: M\to M'$ is CR if
$f_*(HM)\subset HM'$, and $f_*(Jw)=J'f_*(w)$ for every $w\in HM$. A CR
embedding of an abstract CR manifold $M$ into a complex manifold $X$
is a CR map which is an embedding. We say that the embedding is {\it
  generic} if the complex dimension of $X$ is $(m+d)$. 

If $M$ is a real-analytic CR manifold of type $(m,d)$, then by
\cite{ah}, $M$ is locally CR embeddable. Furthermore, by \cite{af}
there exists a complex manifold $X$ such that $M$ can be globally
generically embedded into $X$. Consider a connected open set $\omega$
on $M$. When $M$ is real-analytic, and $f$ is a real-analytic CR
function in $\omega$, then there is a connected open set $\Omega_f$ in
$X$, with $\omega=M\cap \Omega_f$, and a holomorphic extension $\tilde
f$ of $f$ to $\Omega_f$. When $M$ is a generic $C^\infty$-smooth
pseudoconcave CR submanifold of $X$, then there exists a connected
open set $\Omega$ in $X$, with $\omega=M\cap \Omega$, such that any CR
distribution in $\omega$ has a unique holomorphic extension to
$\Omega$. See \cite{bp}, \cite{nv}, \cite{hn3}, \cite{hn4}.

Most of our considerations of
real-analytic CR manifolds will be local, and therefore by the above
mentioned results we can assume without loss of generality that $M$ is
a generically embedded CR submanifold of some open set in $\bb C^n$,
where $n=m+d$. Note here that a compact pseudoconcave CR manifold 
cannot be embedded (even non-generically) into any Stein
manifold (see e.g. \cite{hn}). Therefore, for application purposes the
main theorem is formulated for the target submanifold embedded into
$\bb P^N$. 

Let $M$ be a generic smooth real-analytic submanifold of $\bb C^n$
and let $p\in M$. Then in a sufficiently small neighborhood $U$ of
$p$, $M$ is given by   
\begin{equation}
  M= \{z \in U : \rho_j (z,\overline z)=0,\ \ j=1,\dots d \}, 
\end{equation}
where each $\rho_j$ is a real-valued real-analytic function and
\begin{equation}\label{diff}
  \overline \d\rho_1 \wedge \dots \wedge \overline\d \rho_d \ne 0 {\rm
  \  on \ } M\cap U. 
\end{equation}
We set $\rho=(\rho_1, \rho_2,..., \rho_d)$. There exists a biholomorphic
change of coordinates near $p$, $z=(\xi,\z)\in\bb C^m\times\bb C^d=\bb
C^n$ such  that in the new coordinates, $p=0$, and $M$ is given by
\begin{equation}\label{stdform}
  {\rm Im\  \z}= \phi (\xi,\overline \xi, {\rm Re\ \z}),
\end{equation}
where $\phi$ is a vector-valued real-analytic function with
$\phi(0)=0$ and $d\phi(0)=0$. 

If $U$ is sufficiently small, to every point $w\in U$ we can associate
to $M$ its so-called Segre variety in $U$ defined as 
\begin{equation}\label{Q_w}
  Q_w=\left\{ z\in U: \rho_j(z,\overline w)=0, \ j=1,\dots,d \right\},
\end{equation}
where $\rho_j(z,\overline w)$ is the complexification of the defining
functions of $M$. Another important variety associated with the
submanifold $M$ and the neighborhood $U$ is  
\begin{equation}\label{I_w}
  I_w=\left\{ z\in U: Q_w = Q_z \right\}.
\end{equation}
From the reality condition on the defining functions the following
simple but important properties of Segre varieties follow:
\begin{equation}\label{segre1}
  z\in Q_w \ \Leftrightarrow \ w\in Q_z,
\end{equation}
\begin{equation}\label{segre2}
  z\in Q_z \ \Leftrightarrow \ z\in M,
\end{equation}
\begin{equation}\label{segre3}
  w\in M \Leftrightarrow I_w \subset M.
\end{equation}

If $0\in M$, then from \eqref{diff} and the implicit mapping theorem,
there exist a 
local change of coordinates near the origin, and a pair of small
neighborhoods $U$ and $U_0$ of the origin, $U\Subset U_0$, where $U_0$
is given in the product form 
\begin{equation}\label{standard}
U_0 = {}'U_0\times \tilde U_0,\ \ {}'U_0\subset \bb C^m, \ \ 
\tilde U_0\subset \bb C^d,
\end{equation}
such that for every $w\in U$, the set $Q_w\cap U_0$ can be represented as
the graph of a holomorphic mapping. That is 
\begin{equation}\label{sgraph}
  Q_{w}\cap U_0 = \{ z=(\xi,\z)\in {}'U_0\times \tilde U_0:\
  \z = h(\xi, \overline w)\},
\end{equation} 
where $h(\xi,\overline w)$ is holomorphic in $\xi$ and $\overline
w$. Thus $Q_w$ is a complex submanifold of $U$ of complex codimension
$d$.  
  
The main use of Segre varieties comes from the fact that they are
invariant under biholomorphic mappings. More precisely, given a
holomorphic map $f: U\to U'$, sending a generic smooth real-analytic
submanifold $M$ to another such submanifold $M'$, $f(p)=p'$ implies
$f(Q_p\cap U)\subset Q'_{p'}$ for $p$ sufficiently close to
$M$. An analogous property holds also for holomorphic correspondences.

The proof of the basic properties of Segre varieties in higher 
codimensions is similar to the hypersurface case and can be found in
\cite{ber1} or \cite{me1}.  

A real analytic submanifold $M$ is called essentially finite at $p\in
M$, if $I_{p}=\{p\}$ in a small neighborhood of $p$. The {\it Segre map}
is defined by $\l: w\to Q_w$. A manifold $M$ being essentially finite 
now means that the Segre map is finite near $M$.  It can be shown
(see e.g. \cite{me1}) that any generic Levi non-degenerate CR
submanifold of $\bb C^n$ is essentially finite. Moreover, if $M$ is a 
compact generic submanifold of $\bb C^n$, then it is automatically
essentially finite, since by \cite{df1}, any compact real-analytic
subset of $\bb C^n$ does not contain any non-trivial germs of
complex-analytic varieties. 

We say that $\lambda$ is locally injective at a point $q\in M$, if
there exists a small neighborhood $U_q$ of $q$ such that $\lambda$ is 
an injective map from $U_q$ onto its image. It is easy to see that 
for any Levi non-degenerate hypersurface in $\bb C^n$ the Segre 
map is locally injective.

In \cite{ber} the authors introduced so-called Segre sets. We briefly
recall this construction here. Let $M$ be a generic smooth
real-analytic submanifold of $\bb C^n$, $0\in M$ and let $Q_0 = Q^1_0$
be the usual Segre variety of $0$ as defined in \eqref{Q_w}. Define 
\begin{equation}\label{Q^j}
  Q^j_0 = \bigcup_{z\in Q^{j-1}_0} Q_z, \ \ j>1.
\end{equation}
Then
\begin{equation}
  Q^1_0\subset Q^2_0\subset \cdots \subset Q^j_0\subset\cdots.
\end{equation}
Indeed, let $k$ be the smallest integer such that $Q^k_0\not\subset
Q^{k+1}_0$. Clearly, $k\ge 2$. If $z\in Q^{k}_0\setminus Q^{k+1}_0$,
then there exists $w\in Q^{k-1}_0$ such that
$z\in Q_w$. By assumption, $Q^{k-1}_0\subset Q^{k}_0$. Therefore
$w\in Q^{k}_0$, and $Q_w\subset Q^{k+1}_0$; in particular $z\in
Q^{k+1}_0$, which is a contradiction. 

According to \cite{ber} (see also \cite{ber1} and \cite{mex} for a short
proof of this fact), there exists an
integer $j_0$, $0<j_0<\infty$, such that $\cup_{j\le j_0} Q^j_0$ contains
a neighborhood of the origin in $\bb C^n$, provided that $M$ is
minimal at $0$. We define  
\begin{equation}\label{Omega}
  \O_0 = \{z: z\in Q^{j}_0, \ j\le j_0\}.
\end{equation}
Moreover, if $M$ is compact, or is a relatively compact open set of a
bigger minimal submanifold, then there exists $\e>0$ such that for any
point $p\in 
M$, the neighborhood $\O_p$, defined as in \eqref{Omega}, contains a
ball of radius $\e$ centered at $p$. 

Suppose now that the manifold $M\subset \bb P^n$ is connected and
defined by real polynomials. Then the Segre varieties associated with
$M$ can be defined globally as projective algebraic varieties in $\bb 
P^n$. Indeed, let $M\cap \bb C^n$ be given as a connected component of
the set defined by
\begin{equation}\label{polys1}
  \{z\in \bb C^n : \rho_j(z,\overline z)=0,\ j=1,\dots,d\},
\end{equation}
where $\rho_j$ are real polynomials. We can projectivize each $\rho_j$
to define $M$ in $\bb P^n$ in homogeneous coordinates
\begin{equation}\label{homcoord}
\hat z=[\hat z_0:\hat z_1:\dots:\hat z_n],\ \ 
z_k=\frac{\hat z_k}{\hat z_0},\ \ k=1,\dots,n, 
\end{equation} 
as a connected component of the set defined by
\begin{equation}
  \{\hat z\in \bb P^n : \hat \rho_j (\hat z,\overline{\hat z})=0 \}.
\end{equation}
We may define now the {\it polar} of $M$ as  
\begin{equation}\label{polys2}
  \hat M^c=\{(\hat z,\hat \z)\in \bb P^n \times \bb P^n : \hat
  \rho_j(\hat z,\hat \z)=0,\ \ 
  j=1,\dots,d\}.
\end{equation}
Then $\hat M^c$ is a complex algebraic variety in $\bb P^n\times
\bb P^n$. Given $\t\in \bb P^n$, we set  
\begin{equation}\label{polys3}
  \hat Q_\t = \hat M^c \cap \{(\hat z,\hat\z)\in \bb P^{n}\times\bb P^n
  : \hat\z =\overline \t\}.
\end{equation}
We define the projection of $\hat Q_\t$ to the first coordinate to be
the Segre variety of $\tau$. 

Alternatively, given \eqref{polys1} one can define the polar as
\begin{equation}
  M^c=\{(z,\z)\in \bb C^n \times \bb C^n : \rho_j(z,\z)=0,\ \   
  j=1,\dots,d\}.
\end{equation}
The submanifold $M\cap \bb C^n$ can be recovered by intersecting $M^c$
with the totally real subspace $\mathcal T=\{(z,\z)\in \bb C^{2n} : \z
=\overline z\}$ and taking an appropriate connected component. Given
$w\in \bb C^n$, we define 
\begin{equation}
  Q^c_w = M^c \cap \{(z,\z)\in \bb C^{2n} : \z =\overline w\}.
\end{equation}
The standard Segre variety $Q_w$ can now be recovered by projecting
$Q^c_w$ to $\bb C^n_z$. The algebraic varieties $M^c$ and $Q^c_w$ can
be projectivized which gives us objects geometrically 
equivalent to \eqref{polys2} and \eqref{polys3}. Note that the closure
$\overline{\mathcal T}\subset \bb P^{2n}$ of the set $\mathcal T$ is  a smooth
submanifold of $\bb P^{2n}$, and thus $M$ can be identified with a
connected component of $\overline{M^c}\cap \overline{\mathcal T}$.  

We note that condition \eqref{diff} implies that for $w$ close to $M$,
the Segre variety $Q_w$ contains a connected component $\tilde Q_w$ of
dimension $m$. However in general, $Q_w$ may have other components,
which apriori may even have different dimension (higher than $m$). For
$w\in M$, $w\in \tilde Q_w$, and for $w$ close to $M$, $\tilde Q_w$ is
the component which near $w$ is given by \eqref{sgraph}.  

We will understand essential finiteness of a real algebraic
submanifold $M\subset \bb P^n$ in the sense that for every $w\in M$,
the set  
\begin{equation}\label{I}
\tilde I_w := \{ z\in M:\tilde Q_z =\tilde Q_w\}
\end{equation}
is finite. In this case we can further show that generically the
various $\tilde I_w$ have the same number of points. Note that the set
$\tilde I_w$ is a globally defined object in $\bb P^n$.

\begin{lemma}\label{l-finite}
Let $M$ be a compact smooth real-algebraic essentially finite
generic CR submanifold in $\bb P^n$. Then there exist a neighborhood
$U\subset \bb P^n$ of $M$ and an integer $R\ge 1$ such that 
\begin{equation}
  \#\tilde I_w = \# \{ z\in U: \tilde Q_z =\tilde Q_w\}=R
\end{equation}
for almost all $w\in U$.
\end{lemma} 

\begin{proof} For each $p\in M$ choose a product neighborhood $U_p$ as
  in \eqref{standard}. Let $U=\bigcup_{p\in M} U_p$. We define the map
  $\l$ from $U$ 
  to the set of all algebraic varieties in $\bb P^n$ of dimension
  $m$ by letting $\l(w):=\tilde Q_w$. Then $\tilde Q_w$ depends
  anti-holomorphically on $w$. Locally, near almost every $w\in U$,
  the varieties $\tilde Q_w$ have the same algebraic degree. Since
  algebraic varieties of positive codimension do not divide $U$, the
  $\tilde Q_w$ have fixed degree for almost all $w$ in $U$. Algebraic
  varieties of fixed dimension and degree are known to be parametrized
  by the so-called Chow variety (see e.g. \cite{ha}), and the
  parametrization $\l$ is an algebraic map between algebraic varieties.
  
  Let $Y=\l(U)$. Then since $M$ is essentially finite, $\dim Y=n$. It
  follows (see e.g. \cite{mu}) that there exists an algebraic  
  variety $Z\subset Y$ such that for any $q\in Y\setminus Z$,  
  \begin{equation}\label{deg-fin}  
    \#\l^{-1}(q)=\deg(\l):=R, 
  \end{equation}
  where $R$ is a positive integer. From \eqref{deg-fin} the assertion
  follows. 
\end{proof}

Let $X$ and $X'$ be complex manifolds and $D\subset X$, $D'\subset X'$
be open sets. Recall that if $A\subset D\times D'$ is a holomorphic
correspondence, then $\pi: A\to D$ is proper. $A$ is called a proper
holomorphic correspondence, if $\pi':A\to D'$ is also proper. 
Throughout the paper we will identify the multiple valued map
$F:=\pi'\circ \pi^{-1}$ with its graph $A$, so given a point $p\in D$,  
\begin{equation}\label{cs}
F(p) = \{p'\in D : \pi'\circ \pi^{-1}(p)\}
\end{equation}
is a compact subset of $D'$. Given a complex analytic subset $G'\subset
D'$, $F^{-1}(G')$ is a complex analytic subset of $D$. Indeed
${\pi'}^{-1}(G')$ is clearly analytic, and since $\pi$ is proper, it
follows from Remmert's theorem that $\pi({\pi'}^{-1}(G))$ is analytic
in $D$. $F$ is called a {\it finite-valued} holomorphic
correspondence if $F(p)$ is a finite set for any $p\in D$, and a {\it
  finite} correspondence if in addition $F^{-1}(p')$ is finite for any
$p'\in X'$. If $X'$ is a Stein manifold, then any proper holomorphic
correspondence $F$ is automatically finite-valued. To see this observe
that by the proper embedding theorem, $X'$ can be viewed as a
submanifold of $\bb C^{n'}$ for some $n'>1$, and a compact complex
analytic set $F(p)$ must be discrete. Similarly, if $X$ and $X'$ are
both Stein, then any proper holomorphic correspondence is finite.

Given a finite-valued holomorphic correspondence $A\subset D\times
D'$, there exists 
a complex subvariety $S\subset D$ (possibly empty) such that for any
point $p\in D\setminus S$, there exists a neighborhood $U_p\subset
D\setminus S$, such that $F$ {\it splits} into $k$ holomorphic maps
$F^j: U_p\to D'$, $j=1,\dots,k$, that represent $F$. The integer $k$ is 
independent of $p$, and the $F^j$ are called the {\it branches} of $F$.

Given a locally complex analytic set $A$ in $X$ of pure
dimension $p$, we say that $A$ {\it extends analytically} to an open
set  $U\subset X$, $A\cap U\ne\varnothing$, if there
exists a (closed) complex-analytic set $A^*$ in $U$ such that (i)
$\dim A^*\equiv p$, (ii) $A\cap U\subset A^*$ and (iii) every irreducible
component of $A^*$ has a nonempty intersection with $A$ of dimension
$p$. By the uniqueness theorem for analytic sets such analytic
continuation of $A$ is uniquely defined. From this we define the
analytic continuation of holomorphic correspondences: 

\begin{definition}
Let $D\subset X$ and $ D'\subset X'$ be open sets and let $A\subset
D\times D'$ be a holomorphic correspondence. We say that $A$ extends
as a holomorphic correspondence to an open set $U\subset X$, $D\cap
U\not=\varnothing$, if there exists an open set $U'\subset X'$ such
that $A\cap (U\times U')\ne \varnothing$,  $A$ extends 
analytically to a set $A^*\subset U\times U'$, and $\pi:A^*\to U$ is
proper.   
\end{definition}

If $A$ and $A^*$ are both finite-valued, then $A^*$ may have more 
branches in $D\cap U$ than $A$. The following lemma gives a simple
criterion for the extension to have the same number of branches.

\begin{lemma}\label{l-intersection}
Let $A^*\subset U\times X'$ be a finite-valued holomorphic
correspondence which is an analytic extension of a finite
correspondence $A\subset D\times D'$. Suppose that for any $z\in
(D\cap U)$,  
\begin{equation}\label{pre}
\#\{\pi^{-1}(z)\}=\#\{\pi^{*-1}(z)\},
\end{equation}
where $\pi:A\to D$ and $\pi^*:A^*\to U$ are the projections. Then
$A\cup A^*$ is a holomorphic correspondence in $(D\cup U)\times X'$.
\end{lemma}

The proof is the same as in \cite{sh}, Lemma 2.

\section{Local continuation as a correspondence.}

In this section we prove Theorem \ref{t0}. We may assume the following 
situation: $\O\subset \bb C^n$ is a connected open set, $\O\cap
M=\omega$, $f:\O\to \bb P^N$ is a holomorphic map, and
$f(\O\cap M)\subset M'$. Let $q\in \d \O\cap M$ and let $U$ be 
a neighborhood of $q$  such that for any $z\in U$, the Segre variety
$Q_z$ is defined in a strictly larger set, and can be represented as
in \eqref{sgraph}. We show that $f$ extends to some neighborhood of
$q$ as a holomorphic correspondence. We choose a point $a\in \omega$
such that $df|_{H_a M}$ is an isomorphism, and so close to
$q$ that $q\in \O_{a}$, where $\O_{a}$ is 
defined as in \eqref{Omega}. Fix some neighborhood $U_a\subset \Omega$
of $a$. We will show that $f|_{U_a}$ extends as a holomorphic
correspondence to $\O_a$, in particular to a neighborhood of $q$.      

Assume for simplicity that $a=0$. Choose a small neighborhood $U_0$ of
the origin and shrink $U_a$ in such a way, that for any $w\in U_0$, the set
$Q_w\cap U_a$ is non-empty and connected. Let $U'\subset \bb 
P^N$ be the neighborhood of $M'$ as in Lemma \ref{l-finite}. Define 
\begin{equation}\label{A_0}
  A_0 = \{(w,w')\in U_0\times U' : f(Q_w\cap U_a)\subset Q'_{w'} \}.
\end{equation}
Then $A_0$ is a complex-analytic subset of $U_0\times U'$. In the case
when $M$ and $M'$ are hypersurfaces, this was shown in \cite{sh}, the
proof in the general case is analogous (see also similar constructions
later in this section). Note that in \eqref{A_0},
\begin{equation}\label{eqdim}
\dim Q_w =\dim Q'_{w'}=\dim f(Q_w\cap U_a)=m. 
\end{equation}
Indeed, since $df|_{H_0 M}$ is an isomorphism, $\dim Q_0 =\dim f(Q_0\cap
U_0)$. Without loss of generality we may assume that $M$ near $0$ and
$M'$ near $f(0)$ are chosen as in \eqref{stdform}. Let $J_f(z)$ be the
Jacobian matrix associated with the map $f$. Then in the chosen
coordinate systems, $df|_{H_0 M}$ being an isomorphism means that the
principal minor of $J_f(0)$ of size $m\times m$ has a non-zero
determinant. The same property also holds for points sufficiently
close to the origin. Therefore, after shrinking $U_0$, \eqref{eqdim}
holds for all $w$ in $U_0$. After further shrinking $U_0$ if
necessary, we may assume that if $(w,w')\in A_0$, then $w'\in \tilde
I'_{f(w)}$, which is a finite set by the assumption on $M'$. It
follows that $\dim 
A_0\equiv n$. Finally, \eqref{segre3} implies that if $U_0$ is
sufficiently small, then $A_0$ has no limit points on $U_0\times \d
U'$, and therefore the natural projection from $A_0$ to the first
coordinate is proper. This shows that $A_0$ defines a holomorphic
correspondence. Denote the corresponding multiple valued map by $F_0$.

We shrink the neighborhood $U_0$ of the origin, where $F_0$ is defined,
and choose a 'thin' neighborhood $U_1$ of $Q_0\cap U$ such that for
any $w\in U_1$, the set $Q_w\cap U_0$ is non-empty and connected. Note
that from \eqref{segre1}, $0\in Q_w$ for any $w\in Q_0$. Denote now by
$A_0$ the graph of $F_0$ in $U_0\times \bb P^N$. We define the set
$A_1$ as follows:  
\begin{equation}\label{a1}
  A_1 =  \left\{ (w,w')\in U_1\times \bb P^N : F_0(Q_w\cap U_0)\subset
  Q'_{w'} \right\}. 
\end{equation} 
Then $A_1$ is a complex-analytic subset of $U_1\times \bb P^N$. To
verify this assertion we prove the following: 

\medskip

\begin{enumerate}
\item[1.] $A_1\ne\varnothing$. Indeed, by the invariance property of
  Segre varieties, $f(Q_w\cap U_0)\subset Q'_{f(w)}$ for $w$
  sufficiently close to the origin. Then 
  \begin{equation}\label{e-inc}
    F_0(Q_w\cap U_0)\subset Q'_{w'}, 
  \end{equation}
  where $w'\in F_0(w)$. To see this, suppose that $z\in
  Q_w\cap U_0$ and $z'\in F_0(z)$. By construction, $F_0(Q_z\cap
  U_0)\subset Q'_{z'}$. If $w$ is sufficiently close to the origin,
  $w\in Q_z\cap U_0$, and therefore $w'\in F_0(w)\subset
  Q'_{z'}$. This implies $z'\in Q'_{w'}$. And since 
  $z\in Q_w\cap U_0$ was arbitrary, \eqref{e-inc} holds. Moreover,
  $F_0(w)=\tilde I'_{w'}$, and thus    
  \begin{equation}
    A_0|_{(U_1\cap U_0)\times \bb P^N}\subset A_1,
  \end{equation}
  in particular $A_1$ is non-empty.
  
\item[2.] $A_1$ is a complex-analytic set in a neighborhood of any of
  its points. Indeed, let $(w_0,{w'}_0)\in A_1$, $z_0\in Q_{w_0}\cap
  U_0$ be an arbitrary point, and let $U_{z_0}\subset U_0$ be a small 
  neighborhood of $z_0$. Since $Q_w\cap U_0$ is connected for all
  $w\in U_1$, $F_0(Q_w\cap U_{z_0})\subset Q'_{w'}$  implies
  $F_0(Q_w\cap  U_0)\subset Q'_{w'}$.  Therefore in \eqref{a1} $U_0$
  can be replaced with $U_{z_0}$.
  Choose $U_{w_0}$ so small that $Q_{w}\cap U_{z_0}\ne \varnothing$
  for all $w\in U_{w_0}$. Since the pre-image of an analytic set under
  a holomorphic correspondence is an analytic set,
  $S_{w'}:=F^{-1}_0(Q'_{w'})$ is an analytic subset of $U_0$. Let
  $U'_{w'_{0}}$ be so small that $S_{w'}\cap U_{z_0}\ne \varnothing$
  for all $w'\in U'_{w'_{0}}$. Let $S_{w'}$ near $z_0$ be given by
  \begin{equation}
   S_{w'}=\{z\in U_{z_0}: \phi_j(z,\overline w')=0,\
   j=1,\dots,\tilde j \}, 
  \end{equation}
  where the $\phi_j$ depend holomorphically on $\overline w'$. We may
  assume that $U_{z_0}$ is chosen as in \eqref{standard}, and
  therefore $z\in Q_w$ simply means $z=(\xi, \zeta)$ and
  $\zeta=h(\xi,\overline w)$. Then the condition $F_0(Q_w\cap
  U_{z_0})\subset Q'_{w'}$ is equivalent to  
  \begin{equation}\label{ajsys}
    \phi_j\left((\xi,h(\xi,\overline w)),\overline w'\right)=0, \
    \ \xi\in {}'U_{z_0},\ j=1,\dots,\tilde j.
  \end{equation}
  This is an infinite system of holomorphic equations (after
  conjugation) which defines $A_1$ as an analytic set in
  $U_{w_0}\times U'_{w'_0}$.  
  
\item[3.] $A_1$ is closed in $U_1\times \bb P^N$. Indeed,
  suppose that $(w^j,{w'}^j)\to(w^0, {w'}^0)$, as $j\to\infty$,
  where $(w^j, {w'}^j)\in A_1$ and $(w^0,{w'}^0)\in U_1\times{\bb
  P}^N$. Since $Q_{w^j}\to Q_{w^0}$, and $Q'_{{w'}^j}\to Q'_{{w'}^0}$
  as $j\to\infty$, by analyticity $F_0(Q_{w^0}\cap U_0)\subset
  Q'_{{w'}^0}$, which implies that $(w^0,{w'}^0)\in A_1$ and thus
  $A_1$ is a closed set.   
\end{enumerate}

\medskip

It follows from 1--3 that $A_1$ is a complex-analytic subset of
$U_1\times \bb P^N$. Let $\pi_1:A_1\to U_1$ and $\pi'_1:A_1\to \bb
P^N$ be the coordinate projections. Since $\bb P^N$ is compact,
$\pi_1$ is proper. We consider only the irreducible components of $A_1$ 
of dimension $n$ which contain $A_0$. The union of all such components
we denote again by $A_1$. Thus $A_1$ is an analytic continuation of
$A_0$ as a holomorphic correspondence. Denoting $U_0\cap U_1$ again by
$U_0$, we may assume from the uniqueness theorem for analytic sets
that  
\begin{equation}\label{0=1}
A_1|_{U_0\times \bb P^N} = A_0.
\end{equation}
We set $F_1:=\pi'_1\circ\pi^{-1}_1: U_1\to \bb P^N$. 

By construction, if $(w,w')\in A_1$, then for any $(z,z')\in A_1$ such
that $z\in Q_w\cap U_0$, we necessarily have $z'\in Q'_{w'}$. In order
to construct analytic sets $A_j$ which will extend $A_1$, we wish
to conclude the same for all $z$ in $Q_w\cap U_1$. The difficulty
is that in general, $Q_w\cap U_1$ may have more than one connected
component. To prove the assertion we argue by contradiction, and
assume that there exists a point 
$(w^0,{w'}^0)\in A_1$ and $(z,z')\in A_1$, such that $z\in Q_{w^0}\cap
U_1$,  but $z'\notin Q'_{{w'}^0}$. Connect $w^0$ and the origin with a
smooth path $\g$ contained in $U_1$. We may assume that $w^0$ is the
first point on $\g$ for which the desired property does not
hold. Without loss of generality we may also assume that for all
$p\in\g$  between the origin and $w^0$ (excluding $w^0$) there exists
a small neighborhood $U_p$ such that whenever $z\in U_p$, all
components of $Q_z\cap U_1$ are mapped by $F_1$ into the same Segre
variety. For each point 
$p\in\g$ between the origin and $w^0$, we construct the following set:  
\begin{equation}\label{A_p}
  A_p=\left\{(w,w')\in U(Q_p)\times \bb P^N : F_1(Q_w\cap U_p)\subset
  Q'_{w'} \right\},
\end{equation}
where $U_p\subset U_1$ is a neighborhood of $p$ and $U(Q_p)$ is a
neighborhood of $Q_p\cap U$, which are chosen in such a way that 
$U_p$ satisfies the property described above and that for
any $w\in U(Q_p)$, the set $Q_w\cap U_p$ is connected. Repeating the
argument that was used for $A_1$, one can prove that each $A_p$ is a
complex analytic set, which defines a holomorphic correspondence. For
$p=0$, this is just the set $A_1$. Moreover, for any $p$ between the
origin and $w^0$ we have
\begin{equation}\label{43}
  A_1|_{(U(Q_p)\cap U_1)\times \bb P^N} \subset A_p|_{(U(Q_p)\cap U_1) 
  \times \bb P^N}.
\end{equation}
Indeed, suppose that $w\in U(Q_p)\cap U_1$ and $(w,w')\in A_1$. Let
$z\in Q_w\cap U_p$ be an arbitrary point, and $z'\in F_1(z)$. Then
$F_0(Q_z\cap U_0)\subset Q'_{z'}$. From \eqref{0=1} we have
$F_1(Q_z\cap U_0)\subset Q'_{z'}$. By the assumption on $U_p$,
$F_1(Q_z\cap U_1)\subset Q'_{z'}$, in particular $F_1(w)\subset
Q'_{z'}$. Therefore, $w'\in Q'_{z'}$, and $z'\in Q'_{w'}$. Since $z\in
Q_w\cap U_p$ was arbitrary, it follows that $F_1(Q_w\cap U_p)\subset
Q'_{w'}$. But this means that $(w,w')\in A_p$, and thus \eqref{43}
holds. 

For any $p$, $Q_p\cap U$ is a connected set in $U(Q_p)$ and therefore is
mapped by $A_p$ into the same Segre variety. By continuity and from
\eqref{43} we conclude that $F_1(Q_{w_0}\cap U_1)\subset Q'_{w'_0}$,
which contradicts the assumption. Thus for any $(w,w')\in A_1$, if
$(z,z')\in A_1$ and $z\in Q_w\cap U_1$, then $z'\in Q'_{w'}$.   

\bigskip

We now define recursively for $j>1$ the following sets:
\begin{equation}\label{aj}
  A_j=\left \{(w,w')\in U_j\times \bb P^N : F_{j-1}(Q_w\cap
  U_{j-1})\subset Q'_{w'} \right\}. 
\end{equation}
Here the open set $U_j$ is defined as follows. Suppose that the set
$A_{j-1}\subset U_{j-1}\times \bb P^N$ is already defined and 
$U_{j-1}\subset U$ is some connected open set. We let $U_j$ be the set
of points $w$ in $U$ such that $Q_w\cap U_{j-1}\ne\varnothing$. 
Furthermore, after shrinking at each step, if necessary, the sets
$U_k$ for $k<j$, we may assume that $U_{k-1}\subset U_k$ for $1\le
k\le j$. Note that it follows from the construction that $Q^k_0\subset
U_k$ for $1\le k\le j$. 

We claim that for all $j>0$, $A_j$ is a complex-analytic subset of
$U_{j}\times \bb P^N$, which satisfies the following properties:

\medskip

\begin{enumerate}
\item[(i)]  $A_j|_{U_{j-1}\times \bb P^N} = A_{j-1}$; 
  
\item[(ii)] $A_j$ defines a holomorphic correspondence
  $F_j: U_j\to \bb P^N$;
  
\item[(iii)] for any $(w,w')\in A_j$, if $(z,z')\in A_j$ and $z\in
  Q_w\cap U_j$, then $z'\in Q'_{w'}$.
  
\end{enumerate}

\medskip
Condition (iii) can be understood in the sense that the map $F_j$,
associated with $A_j$, sends all connected components of $Q_{w}\cap
U_j$ into $Q'_{w'}$ provided that $(w,w')\in A_j$. 

The proof is by induction. The case $j=1$ is already proved. Suppose
that $A_{j-1}$ is as claimed. We show that the set defined by
\eqref{aj} is also a holomorphic correspondence satisfying properties
(i)--(iii). 

\smallskip

(i) Let $w\in U_{j-1}$, and $(w,w')\in A_{j-1}$. Then by definition,
$F_{j-2}(Q_w\cap U_{j-2})\subset Q'_{w'}$. From property (i), which by
the induction hypothesis holds for $F_{j-1}$, the correspondences
$F_{j-2}$ and $F_{j-1}$ agree in $U_{j-2}$, and therefore we have
\begin{equation}\label{e-f}
F_{j-1}(Q_w\cap U_{j-2})\subset Q'_{w'}.
\end{equation}
From (iii), $F_{j-1}$ maps all components of $Q_w\cap U_{j-1}$ into
the same Segre variety. Therefore \eqref{e-f} implies $F_{j-1}(Q_w\cap 
U_{j-1})\subset Q'_{w'}$, which by definition means that $(w,w')\in
A_j$. In particular, the set $A_j$ is non-empty. Condition (i) for
$A_j$ will be completely proved, once we know that $A_j$ is a
complex-analytic set, $\dim A_j\equiv n$, and select only the
irreducible components of $A_j$ which have intersection with $A_{j-1}$
of dimension $n$. 

\smallskip

Proof of (ii). Let $(w^0,{w'}^0)\in A_j$. If $Q_{w}\cap U_{j-1}$ is  
connected for all $w$ sufficiently close to $w^0$, then the proof of
the fact that $A_j$ is complex-analytic 
near $(w^0,{w'}^0)\in A_j$ is the same as for $A_1$ in Step 2. Let
$\tilde U$ be the largest connected relatively open subset of $U_{j}$
such that $0\in \tilde U$ and
for all $z\in \tilde U$, $F_{j-1}$ maps $Q_z\cap U_{j-1}$ into the
same Segre variety. From property (iii) for $A_{j-1}$ we have
$U_{j-1}\subset \tilde U$. Then \eqref{aj} defines a holomorphic
correspondence $\tilde A\subset\tilde U\times \bb P^N$. The proof is
the same as in Step 2 for $A_1$. Denote by $\tilde F$ the multiple
valued map associated with $\tilde A$. By repeating the argument used
for $A_1$, we can show that for any $w\in \tilde U$, $\tilde F$ maps
all connected components of $Q_w\cap \tilde U$ into the same Segre 
variety. 

For each point $p\in \tilde U$ we may define now the following set
\begin{equation}\label{prep2}
  A_{p} = \{(w,w')\in U(Q_{p})\times \bb P^N :
  \tilde F(Q_w\cap U_{p})\subset Q'_{w'} \},
\end{equation}
where the neighborhoods $U(Q_{p})$ of $Q_p$ and $U_p$ of $p$ are
chosen as in the construction of the set defined by
\eqref{A_p}. Let $F_p$ be the map associated with $A_p$.
Clearly, $F_{p}$ coincides with $F_{j-1}$ for $p$ sufficiently close
to the origin. We claim that the map $F_p$ agrees with $F_{j-1}$ in
$U(Q_p)\cap U_{j-1}$ 
for all $p\in \tilde U$. Indeed, suppose that $w\in U(Q_p)\cap
U_{j-1}$ and $(w,w')\in A_{j-1}$. Then $F_{j-2}(Q_w\cap
U_{j-2})\subset Q'_{w'}$. To prove the assertion we need to show that 
\begin{equation}\label{kill-it}
  \tilde F(Q_w\cap U_p)\subset Q'_{w'}.
\end{equation}
Let $z\in Q_w\cap U_p$ be an arbitrary point, and $z'\in \tilde
F(z)$. Then $F_{j-1}(Q_z\cap U_{j-1})\subset Q'_{z'}$. Since $\tilde
F$ and $F_{j-1}$ agree in $U_{j-1}$, it follows that $\tilde F(Q_z\cap 
U_{j-1})\subset Q'_{z'}$. For $z\in \tilde U$, $\tilde F$ maps
different components of $Q_z\cap \tilde U$ into the same Segre
varieties, and therefore we have $\tilde F(Q_z\cap \tilde U)\subset
Q'_{z'}$. In particular, $w'\in Q'_{z'}$, which implies $z'\in
Q'_{w'}$. Since $z$ was arbitrary, \eqref{kill-it} holds. 

By the construction, $F_p$ maps $Q_p$ into the same Segre variety for
all $p\in \tilde U$. By analyticity this means that for any point $w$
in $\d \tilde U\cap U_{j}$, $F_{j-1}$ maps $Q_w\cap U_{j-1}$ into the
same Segre variety. Therefore, $\tilde U=U_{j}$. We choose only
irreducible components of $A_j$ of dimension $n$ which contain
$A_{j-1}$. This proves (ii) and also completes the proof of (i).  

\smallskip

Finally, property (iii) can be shown the same way as it was done for
$A_1$.  

\bigskip

By the construction, from minimality of $M$ and from \cite{ber}, for
some $j_0>1$, the set $A_{j_0}$ defines a holomorphic correspondence
$F_{j_0}$ in a neighborhood $\Omega_0\subset U_{j_0}$ of the origin. 
Note that the size of this neighborhood depends only on the geometry
of $M$ and is independent of $U_0$, where $f$ was originally
defined. It remains to show now that $F_{j_0}$ satisfies
\begin{equation}\label{last}
F_{j_0}(M\cap \Omega_0)\subset M'.
\end{equation}
If $(z,z')\in A_{j_0}$, then $F_{j_0-1}(Q_z\cap U_{j_0-1})\subset
Q'_{z'}$. From property (i), $F_{j_0}(Q_z\cap U_{j_0-1})\subset
Q'_{z'}$, and from (iii) it follows that 
\begin{equation}\label{prop-incl}
F_{j_0}(Q_z\cap U_{j_0})\subset Q'_{z'}.
\end{equation}
Suppose now that for some $z^0\in M\cap \Omega_0$,
$F_{j_0}(z^0)\not\subset M'$. Then there exists $z'\in
F_{j_0}(z^0)\setminus M'$. Note that $F_{j_0}(M\cap U_0)\subset M'$,
and therefore by continuity we may find $z^0$ and $z'$ such that $z'$
is close to $M'$. From \eqref{prop-incl} we have
$F_{j_0}(Q_{z^0}\cap U_{j_0})\subset Q'_{z'}$. Since $z^0\in Q_{z^0}$,
we have $F_{j_0}(z^0)\subset Q'_{z'}$, in particular, $z'\in
Q'_{z'}$. But from \eqref{segre2}, $z'\not\in Q'_{z'}$, since
$z'\not\in M'$. This contradiction proves \eqref{last}. 

Theorem \ref{t0} is proved. Note that in general, $F_{j_0}$ may not be 
finite-valued. However, by the Cartan-Remmert theorem (see e.g. \cite{lo}) 
combined with Remmert's proper mapping theorem, the set of points 
\begin{equation}
  \S = \{ z\in U_{j_0} : \dim \pi^{-1} (z) >0\}
\end{equation}
is a complex subvariety of $U_{j_0}$. Since $\dim A_{j_0}\equiv n$, 
we have $\dim \S<n$, in particular, $\S$ does not divide $U_{j_0}$, 
and  $A_{j_0}|_{(U_{j_0}\setminus \S)\times \bb P^N}$ is a finite-valued
holomorphic correspondence.

\section{Extension as a finite correspondence along $M$.}

In this section we give the proof of Theorem \ref{t1}. One difficulty
in the proof of analytic continuation of correspondences lies in the
fact that the continued correspondence may acquire additional
branches. To deal with this we define the notion of a {\it complete}
correspondence as follows.

\begin{definition}
  Let $M\subset \mathbb C^n$ be a smooth real-analytic essentially
  finite generic CR submanifold, and let $M'\subset \mathbb P^N$ be a
  smooth compact  real-algebraic essentially finite
  generic submanifold. Let $F: U\to \mathbb P^N$ be a holomorphic
  correspondence such that $F(U\cap M)\subset M'$. Then $F$ is
  called {\sl complete} if for every $z\in M$, we have
  $F(z)={\tilde I'}_{z'}$ for some (and therefore for any) $z'\in
  F(z)$.  Here $\tilde I'_{z'}$ is defined as in \eqref{I}.
\end{definition}

Note that since $M'$ is essentially finite, a complete correspondence
is finite-valued near $M$, but in general it may be reducible, even if
defined on all of $M$.

Assume that $M$ is locally, near $p\in M$, generically embedded into
an open set in $\bb C^n$, so $n=m+d$. Let $f$ be a holomorphic map
defined in a neighborhood $U_p\subset \bb C^n$ of $p\in M$, of maximal
rank and such that $f(U_p\cap  M)\subset M'$. We replace $f$ with a
complete correspondence. For that we choose a small neighborhood 
$U_0$ of $p$ and shrink $U_p$ in such a way, that for any $w\in U_0$,
the set $Q_w\cap U_p$ is non-empty and connected. We define $A_0$ as
in \eqref{A_0}. Denote the corresponding multiple valued map by $F_0$.
Then $F_0$ is a complete holomorphic correspondence. Indeed, since $f$
is of maximal rank and $d'\ge d$, for $z\in U_0$, $\dim f(Q_z\cap
U_p)=\dim Q_z=m$. For $z$ sufficiently close to $p$, $f(z)$ is close
to $M'$, and since $M'$ is essentially finite, the only points whose
Segre varieties can contain $f(Q_z\cap U_p)$ are in
$I'_{f(z)}$. Furthermore, $F_0$ is a finite 
correspondence. To see this we observe that if $E\subset U_0$ is a
positive dimensional set such that $F_0(E)$ is discrete, then by the
construction $f(Q_z\cap U_p)\subset Q'_{z'}$ for all $z\in E$. But
since $\bigcup_{z\in E} Q_z\cap U_p$ has dimension bigger than $m$,
this contradicts the fact that $f$ is of maximal rank in $U_p$. Therefore
the pre-image of any point under $F_0$ is finite, and thus $F_0$ is a
finite correspondence. 

We now show that $F_0$ extends as a finite correspondence along any
path on $M_1$. Our construction of the analytic continuation of $F_0$
will preserve completeness, and therefore from Lemmas \ref{l-finite} 
and \ref{l-intersection} we conclude
that such analytic continuation will have the same number of branches
near $M_1$.  The problem can be localized as follows. Let $\g: [0,1]\to
M_1$ be the given path, $\g(0)=p$, and assume that for $t_0\le 1$,
$q=\g(t_0)$ is the first point on $\g$ to which $F_0$ does not
extend as a finite correspondence. We choose $t_1\in [0,t_0)$ so
close to $t_0$, that for $a=\g(t_1)$ we have $q\in \Omega_{a}$, where
$\Omega_a$ is defined as in \eqref{Omega}. Then by Theorem \ref{t0},
$F_0$ extends as a holomorphic correspondence $A\subset \Omega_a\times
\bb P^n$. Thus we only need to show that after possibly shrinking 
the neighborhood $\Omega_a$, the set $A$ is a finite correspondence. 
Let $\pi: A\to \Omega_a$ and $\pi': A\to \bb P^N$ be the natural
projections, and set $F:=\pi'\circ\pi^{-1}$.   

From \eqref{last} there exists a neighborhood $\tilde U\subset
\Omega_a$ of $M\cap \Omega_a$ such that $F(\tilde U)\subset U'$, where
$U'\subset \bb P^N$ is a neighborhood of $M'$ as in Lemma
\ref{l-finite}. We now 
repeat the argument of analytic continuation of $A_0$ along $\gamma$
by constructing the sets $A^*_j$, $j=1,2,\dots$, with the only
difference that the standard neighborhoods of the form
\eqref{standard} are chosen so small that they are contained in
$\tilde U$. Since the new set $\O_{a}$ may be smaller than $\tilde U$,
this continuation may require more than one step. More precisely, we
choose a sequence of points $\{a_\nu\}_{\nu=0}^l$ such that $a_\nu \in
\gamma$, $a_0=a$, $a_l=q$ and $a_{\nu}\in \O_{a_{\nu-1}}$ for
$0<\nu\le l$. For each $a_\nu$ starting with $a_0$ we use Theorem
\ref{t0} to extend a finite correspondence $F_\nu$ defined in a
neighborhood of $a_\nu$ to a holomorphic correspondence $F_{\nu+1}$
defined in $\Omega_{a_\nu}\subset \tilde U$. This time we show in
addition that at each step the extension is a finite
correspondence. Since $a_{\nu+1}\in \O_{a_\nu}$, the process can be
continued until we reach the point $q$.  

Suppose that $F_{\nu+1}:\Omega_{a_\nu}\to \bb P^N$ is a holomorphic
correspondence, which is obtained from the finite correspondence $F_{\nu}$
defined in a small neighborhood of the point $a_\nu$ for some $\nu$ by
the inductive construction of the sets $A^*_j$; that is
\begin{equation}
  A^*_j=\left\{(w,w')\in U^*_j\times \bb P^N: F^*_{j-1}(Q_w \cap
  U^*_{j-1})\subset Q'_{w'}\right\},
\end{equation}
where the Segre set $Q^j_{a_\nu}$ is contained in $U^*_j$, the map
$F_{\nu}$ is associated with the set $A^*_0$, and $F_{\nu+1}$
corresponds to the set $A^*_{j_0}$ with $U^*_{j_0}=\Omega_{a_\nu}$. 
Let $F^*_j$ be the map associated with $A^*_j$. 

Clearly, the sets $A^*_j$ are contained in the set $A$, which defines
a correspondence $F:\tilde U\to U'$. We claim that $F_{\nu+1}$ is a
finite correspondence. To prove this assertion we let $k$ be the
smallest integer such that $F^*_k$ is not a finite correspondence. By
assumption, $k>0$. Suppose that there exists a point $w'\in U'$ such
that the analytic set ${F^*_k}^{-1}(w')\subset U^*_k$ has positive
dimension. By the construction we have
\begin{equation}\label{bad}
  F^*_{k-1}(Q_z\cap U^*_{k-1})\subset Q'_{w'}, 
  {\rm \ \ for\ all\ \  } z\in {F^*_k}^{-1}(w').
\end{equation}
Since $M$ is essentially finite, there exists a set $E\subset
{F^*_k}^{-1}(w')$ such that $\dim \cup_{z\in E} Q_z=m+1$. It follows
from \eqref{bad} that $\dim F^*_{k-1}(\cup_{z\in E} Q_z) \le m$. But
this contradicts  the induction hypothesis that $F^*_{k-1}$ is
finite. 

Suppose now that there exists a point $w\in U^*_k$ such that
$F^*_k(w)$ is not discrete. By the construction this means that
$F^*_{k-1}(Q_w\cap U^*_{k-1})\subset Q'_{z'}$, where $z'$ belongs to
the non-discrete set. Since $F^*_{k-1}$ is finite, $F^*_{k-1}(Q_w\cap
U^*_{k-1})$ has dimension $m$. But then there can be only finitely
many $z'$ whose Segre varieties can contain $F^*_{k-1}(Q_w\cap
U^*_{k-1})$. Therefore $F^*_k$ is also finite for any $k\ge 0$.

Hence, the obtained extension $F_{\nu+1}$ is a finite
correspondence. We can repeat the same argument for extending $F_0$
along $\gamma\cap \tilde U$ until we reach the point $q$.

Thus we have proved that $F_0$ extends along any path on $M_1$ as a
finite correspondence $F$. It follows from the construction that
the number of branches of $F$ coincides with the number $R$ defined in
Lemma \ref{l-finite} for almost all points on $M_1$. It remains to
observe now that from the simple connectivity of $M_1$ and the
monodromy theorem, the extension of $F_0$ along homotopically
equivalent paths gives the same result. The version of the monodromy
theorem for finite-valued correspondences can  be found in \cite{sh}.    

Part (2) of the statement of the theorem follows from the construction
of the extension $F$. Indeed, from the construction, if $z\in M$, and
$z'\in F(z)$, then $F(z)=\tilde I_{z'}$. If $\lambda'$ is locally
injective near $z'$, then $F$ splits into $R$ holomorphic mappings near 
$z$. 

The proof of Theorem \ref{t1} is now complete. 


\section{Pseudoconcave submanifolds in $\bb P^n$.}

In this section we prove the rest of the results stated in Section 2.

\begin{proof}[Proof of Theorem \ref{t1.5}] $\ $
  \begin{enumerate}
  \item[(a)] According to Thm 5.2 of \cite{hn2}, there exists an $m+d$
    dimensional  irreducible algebraic subvariety $Y$ of $\bb P^n$ such
    that $M$ is a generic CR submanifold of the regular part of $Y$,
    ${\rm reg\ } Y$. Because of the pseudoconcavity of $M$ (or because
    of Property E), the continuous CR map $f$ is smooth and has a
    unique holomorphic extension to an open neighborhood $\O$ of $M$
    in ${\rm reg\ } Y$. Thus $f$ can be regarded as a holomorphic map
    from $\Omega$ to $\bb P^N$. This means that $f$ may be given by
    $N$ meromorphic functions $f_1, f_2, \dots f_N$ in
    $\Omega$. To see this, we choose homogeneous coordinates
    $[z_0:z_1:\dots:z_N]$ in $\bb P^N$ such that the hyperplane
    $\{z_0=0\}$ is in general position with respect to $f(\O)$,  and
    set $\O_j=f^{-1}(V_j)$, where $V_j=\bb P^N\setminus\{z_j=0\}$,
    $j=0,1,\dots,N$. The $\O_j$ give an open cover of $\O$, in each
    $V_j$ we have the inhomogeneous coordinates
    $(w_{1j},w_{2j},\dots,w_{Nj})$, where $w_{ij}=z_i/z_j$, and
    $f|_{\O_j}$ is given by holomorphic functions 
    \begin{equation}
      \Omega_j\ni t \to (w_{1j}(t),w_{2j}(t), \dots, w_{Nj}(t)).
    \end{equation}
    We define the meromorphic functions $f_1, f_2,\dots, f_N$ by
    $f_k(t)=w_{k0}(t)$ in $\Omega_0$, and by
    $f_k(t)=w_{kj}(t)/w_{0j}(t)$ in $\Omega_j$ for
    $j,k=1,2,\dots,N$. Note that these definitions are consistent on
    the overlaps.

    By Thm 5.2 of \cite{hn2}, each $f_k$ is the restriction to $M$
    of a rational function on $Y$, and hence can be regarded as a
    rational function on $\bb P^n$. This gives the desired rational
    map from $\bb P^n$ to $\bb P^N$. 
    
  \bigskip

  \item[(b)] Since $M$ is generic in $\bb P^n$, and since $f$ is a
  local CR diffeomorphism, the Jacobian $\det
  J_f$ of the extension of $f$ to $\bb P^n$ is not identically zero.
  Hence, the set $\S = \{z\in \bb P^n: \det J_f(z)=0\}$ if non-empty, is
  a subvariety of $\bb P^n\setminus \L$ of complex codimension one,
  where $\L$ is the indeterminacy locus of $f$. Suppose $\S\ne
  \varnothing$. Then, since $f$ is locally biholomorphic near any
  point on $M$, $M\cap \S=\varnothing$. On the other hand, by the
  Remmert-Stein theorem, the closure of $\S$ is a subvariety of $\bb
  P^n$ of codimension one. It is well known that its complement in
  $\bb P^n$ is therefore a Stein manifold. But a pseudoconcave $M$ (or
  an $M$ satisfying  Property E) has no CR embedding into a Stein
  manifold (see \cite{hn} and \cite{hn2}). Thus $\S=\varnothing$.     
  
  Let $F:\bb C^{n+1}\to \bb C^{n+1}$ be a polynomial map such that
  $f\circ\pi = \pi\circ F$, where $\pi: \bb C^{n+1}\to \bb P^n$ is the
  canonical projection. Without loss of generality assume that the
  components of $F$ are homogeneous polynomials of degree $k$ without
  common factors. We claim that $\det J_F(z)\ne 0$ for any point $z\in
  \bb C^{n+1}$. Indeed, suppose on the contrary that 
  \begin{equation}
    E=\{z\in \bb C^{n+1}: \det J_F(z)=0\}
  \end{equation} 
  is a non-empty subvariety of complex codimension one. Then
  $F(E)\ne\{0\}$, and therefore there exists a point $p\in E$ such
  that $F(p)\ne 0$. For $z\in \bb C^{n+1}$, let $L_z:=\{\l z: \l\in
  \bb C\}$ be the complex line passing through the point $z$ and the
  origin. Since the Jacobian of $f$ does not 
  vanish outside the indeterminacy locus, and $F(L_p)\ne\{0\}$, there
  exists a small neighborhood $U$ of $p$ such that for all $z$ and $w$
  in $U$, $z\ne w$,  
  \begin{equation}\label{inj}
    F(L_z)\cap F(L_w)=\{0\}.  
  \end{equation}
  Furthermore, $F|_{U\cap L_z}=\l^k z$, and after shrinking $U$ if
  necessary, we may assume that $F|_{U\cap L_z}$ is an injective
  function for all $z\in U$. From this and \eqref{inj} we
  conclude that $F$ is injective in U, which contradicts the
  assumption that $p\in E$. Thus $\det J_F\ne 0$ and therefore is a
  constant. 

  Finally, observe that $\det J_F(z)$ is a homogeneous polynomial of
  degree $(k-1)(n+1)$, and being constant means that $k=1$, i.e. $F$
  is a linear automorphism.
  \end{enumerate}
\end{proof}

\begin{proof}[Proof of Theorem \ref{t2} and Corollary \ref{c3}.] 
  We may regard $M$ as being a generic CR submanifold of a complex
  manifold $X$. Note that the pseudoconcavity of $M$ (or Property E)
  implies that $M$ is minimal, because minimality is well-known to be
  equivalent to wedge extendability. There is a neighborhood $V_p$ of
  $p$ in $X$ such that $f$ extends to a holomorphic mapping $f:V_p\to
  \bb P^N$. It is easy to check that, by possibly shrinking $V_p$, the
  extended map $f$ has maximal rank in $V_p$. Hence we may conclude
  from Theorem \ref{t1} that $f$ extends to a finite holomorphic
  correspondence $F: V\to \bb P^N$, where $V\subset X$ is some
  neighborhood of $M$. 

  Since the Segre map associated with $M'$ is injective, $F$ splits at
  every point of $M$, and every map $F^j$ of the splitting is a CR map
  from $M$ to $M'$. One of them, say $F^1$ is the extension of $f$. By
  Theorem \ref{t1.5}(a), $F^1$ extends to a rational map from $\bb P^n$
  to $\bb P^N$. Moreover, if $n=N$ and $M$ is generic in $\bb P^n$,
  then $F^1$ is locally biholomorphic, because the Segre map
  associated with $M$ is injective. Thus by Theorem \ref{t1.5}(b)
  $F^1$ extends to a linear automorphism of $\bb P^n$. Rationality of
  $F^1$ implies that $M$ must also be algebraic,  which proves
  Corollary \ref{c3}. 
\end{proof}

\begin{proof}[Proof of Theorem \ref{c}.]
  Since $M$ and $M'$ are Levi non-degenerate, the associated Segre
  maps are locally injective, and $M$ and $M'$ satisfy the conditions
  of Theorem~\ref{t1}. If $M$ and $M'$ are strictly pseudoconvex then
  $f$ extends as a locally biholomorphic map to a neighborhood of $p$ by
  \cite{pt}. If $M$ is pseudoconcave, the result follows from
  \cite{hn1}. Therefore the map $f$ defined near $p$ also satisfies
  the conditions of Theorem~\ref{t1}. Thus $f$ extends as a finite
  correspondence along $M$. Since the set $\Sigma'$, where the Segre
  map associated with $M'$ branches, is empty, the extended correspondence is
  single-valued. Since the Segre map associated with $M$ is injective,
  the extension $f:M\to M'$ is a locally biholomorphic map.

  We now show that the extension is globally biholomorphic in a
  neighborhood of $M$. For that we note that $M'$ is simply
  connected. Indeed, if $k=n-1$ or $k=0$, then $M'=S^{2n-1}$ which is
  simply connected. If $0<k<n-1$, then we choose an affine patch $V'_n$
  in $\bb P^n$ where $z'_n\ne 0$. Then
  \begin{equation}
    M'\cap V'_n =\{|w'_0|^2 +\dots+ |w'_k|^2-|w'_{k+1}|^2
    -\dots-|w'_{n-1}|^2 =1\}, 
  \end{equation}
  where $w'_j=z'_j/z'_n$. Let $\pi$ be the projection from $M'$ to the
  coordinates $(w'_{k+1},\dots,w'_{n-1})$. Then $\pi$ is onto, and for
  any $w=(w'_{k+1},\dots,w'_{n-1})$, $\pi^{-1}(w) \cong S^{2k+1}$,
  which is simply connected. Therefore, $M'$ is also simply
  connected. Because of that the germ of a map $f^{-1}$ extends
  holomorphically along any path on $M'$ to a holomorphic map $f^{-1}:
  M'\to M$. Thus $f|_{M}$ maps $M$ one-to-one and onto $M'$, and $f$
  is globally biholomorphic.

  If $M'$ is pseudoconcave, then by the invariance of the Levi form,
  so is $M$. Thus by Theorem~\ref{t1.5}(b) $f$ extends to a linear
  automorphism of $\bb P^n$.   
\end{proof}


\end{document}